\documentclass[12pt]{amsart}

\usepackage[scale=0.75]{geometry} 
\usepackage{amsmath}
\usepackage{amssymb}
\usepackage{latexsym}
\usepackage{amsthm}
\usepackage{mathrsfs}
\usepackage{color}
\usepackage{amsfonts}
\usepackage{enumitem}
\usepackage{array}
\usepackage{setspace}
\usepackage{mathtools}
\usepackage{tikz}
\usepackage{tikz}
\usetikzlibrary{shapes.misc,shadows}
\usetikzlibrary{positioning}
\usepackage{braids}
\usepackage{tikz-cd}
\usepackage{hyperref}
\usetikzlibrary{arrows}
\usetikzlibrary{knots,patterns}

\usepackage{stmaryrd,graphicx}

\newtheorem{theorem}{Theorem}[section]
\newtheorem*{theorem*}{Theorem}

\newtheorem*{cor*}{Corollary}

\newtheorem*{lem*}{Lemma}
\newtheorem{prop}[theorem]{Proposition}
\newtheorem*{prop*}{Proposition}

\theoremstyle{definition}
\newtheorem{definition}[theorem]{Definition}
\newtheorem*{definition*}{Definition}
\newtheorem{conjecture}[theorem]{Conjecture}
\newtheorem*{conjecture*}{Conjecture}
\newtheorem*{condition*}{Condition}
\newtheorem*{assumption*}{Assumption}

\theoremstyle{remark}

\newtheorem*{rem*}{Remark}

\newtheorem{problem}[theorem]{Problem}
\newtheorem*{problem*}{Problem}
\newtheorem*{question*}{Question}

\numberwithin{equation}{section}

\newcommand{\BC}{\mathbb C}
\newcommand{\BF}{\mathbb F}

\newcommand{\BZ}{\mathbb Z}
\newcommand{\BP}{\mathbb P}

\newcommand{\CE}{\mathcal E}

\DeclareMathOperator{\Jac}{Jac}

\DeclareMathOperator{\GL}{GL}
\DeclareMathOperator{\SL}{SL}

\DeclareMathOperator{\Hom}{Hom}

\DeclareMathOperator{\res}{res}

\DeclareMathOperator{\pExp}{Exp}
\DeclareMathOperator{\pLog}{Log}
\DeclareMathOperator{\trace}{tr}

\title{Cohomology rings of character varieties}

\author{Anton Mellit}
\address[Anton Mellit]{Faculty of Mathematics, University of Vienna, Vienna, Austria}
\email{anton.mellit@univie.ac.at}

\date{\today}

\begin{document}
	\onehalfspacing
	
	\begin{abstract}
		In this talk I give an introduction and present some recent progress towards understanding the cohomology rings of character varieties of Riemann surfaces, such as the proof of the $P=W$ conjecture and the computation of the zero-dimensional COHA. In the case of punctured sphere I present an explicit description relating the cohomology rings to the Hilbert scheme of $\BC^2$, refining conjectures of Hausel-Letellier-Rodriguez-Villegas and Chuang-Diaconescu-Donagi-Pantev. I explain how the general case should be related to the symplectic geometry of the Hilbert scheme.
	\end{abstract}
	
	\maketitle
	
	
	\section{Classical motivation}
	One way character varieties naturally arise is as follows. Suppose we have a linear ODE
	\begin{equation}\label{eq:ODE}
		\frac\partial{\partial t} \vec{\varphi}(t) = A(t) \vec{\varphi}(t)\quad(t\in \mathbb{C}).
	\end{equation}
	Here $A(t)$ is an $n\times n$ matrix-valued meromorphic function of one complex variable $t$, and we are looking for vector-valued solutions $\vec{\varphi}(t)$. Let $t_1,\ldots, t_k$ be the poles of $A(t)$ (including $\infty$ if necessary), which we call \emph{punctures}. Choose a basepoint $\ast$, away from the poles of $A(t)$. In a neighborhood of $\ast$, we have an $n$-dimensional vector space of solutions of the ODE \eqref{eq:ODE}. For any loop $\gamma$ based at $\ast$ avoiding the punctures, solving the ODE along $\gamma$ produces a linear transformation of the space of solutions. Choosing a basis, we obtain the \emph{monodromy representation} of the fundamental group
	\[
	\rho:\pi_1(*, \mathbb{P}^1(\mathbb{C})\setminus\{t_1,\ldots,t_k\}) \to \mathrm{GL}_n(\mathbb{C}).
	\] 
	A different choice of basis produces an equivalent representation.
	
	Another way to think about this is that local solutions to the ODE form a \emph{locally constant sheaf} on $\mathbb{P}^1(\mathbb{C})\setminus\{t_1,\ldots,t_k\}$. Locally constant sheaves are also called \emph{local systems}, and they are in a natural bijection with equivalence classes of monodromy representations.
	
\begin{figure}\label{fig:gammas}
		\begin{tikzpicture}[scale=1.]
			\draw (2.7,3.5) circle (2.5);
			\draw (3-0.5, 3+2.5) node[draw,circle, inner sep=2pt,fill] (b) {} node[above right] {$*$};
			\draw (3-2, 3-.5) node[draw,circle, inner sep=2pt,fill] (x1) {} node[above] {$t_1$};
			\draw (3-0.5, 3-.9) node[draw,circle, inner sep=2pt,fill] (x2) {} node[above] {$t_2$};
			\draw (3, 3-1.2+0.0) node[above right] {$\cdots$};
			\draw (3+1.3, 3-.6) node[draw,circle, inner sep=2pt,fill] (xk) {} node[above] {$t_k$};
			\coordinate[left=0.15 of x1] (x1L);
			\coordinate[right=0.15 of x1] (x1R);
			\coordinate[left=0.15 of x2] (x2L);
			\coordinate[right=0.15 of x2] (x2R);
			\coordinate[left=0.15 of xk] (xkL);
			\coordinate[right=0.15 of xk] (xkR);
			\draw[>=latex,dashed,->] (b) to [out=-120, in=90]  node[near end, left] {$\gamma_1$} (x1L);
			\draw[>=latex,dashed,->] (x1L) to [out=-90, in=-90, looseness=3] (x1R); 
			\draw[>=latex,dashed] (x1R) to [out=90, in=-110] (b);
			\draw[>=latex,dashed,->] (b) to [out=-90, in=90]  node[near end, left] {$\gamma_2$} (x2L);
			\draw[>=latex,dashed,->] (x2L) to [out=-90, in=-90, looseness=3] (x2R); 
			\draw[>=latex,dashed] (x2R) to [out=90, in=-80] (b);
			\draw[>=latex,dashed,->] (b) to [out=-60, in=90]  node[near end, left] {$\gamma_k$} (xkL);
			\draw[>=latex,dashed,->] (xkL) to [out=-90, in=-90, looseness=3] (xkR); 
			\draw[>=latex,dashed] (xkR) to [out=90, in=-50] (b);
		\end{tikzpicture}
		\caption{Fundamental group of a $k$-punctured sphere}
	\end{figure}
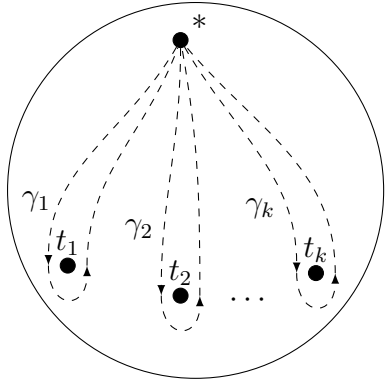

	More explicitly, the fundamental group of $\mathbb{P}^1(\mathbb{C})\setminus\{t_1,\ldots,t_k\}$ is the free group on $k-1$ generators. Drawing loops around the punctures as shown on Figure \ref{fig:gammas}, one obtains the following presentation of the fundamental group:
	\begin{equation}\label{eq:presentation pi genus 0}
		\pi_1(*, \mathbb{P}^1(\mathbb{C})\setminus\{t_1,\ldots,t_k\}) = \langle \gamma_1, \ldots, \gamma_k \;|\; \gamma_1 \cdots \gamma_k = 1 \rangle.
	\end{equation}
	This presentation is better than the free group presentation for the following reason. Let us assume that $A(t)$ has simple poles. In this case, analyzing the ODE locally around each puncture, one can compute the \emph{local monodromy} matrices $\rho(\gamma_i)$, up to conjugation, by a purely algebraic procedure. For instance, the characteristic polynomial of $\rho(\gamma_i)$ coincides with the characteristic polynomial of the matrix $\exp{\left(2\pi \sqrt{-1} A_i\right)}$, where $A_i$ is the residue of $A(t)$ at $t=t_i$.
	
	In view of the presentation \eqref{eq:presentation pi genus 0}, a monodromy representation is encoded by a $k$-tuple of $n\times n$ matrices $M_1, \ldots, M_k$ satisfying
	\[
	M_1 \cdots M_k = I_n,
	\]
	where the local monodromy computation informs us about the conjugacy class of each $M_i$. 
	
	Contrary to the local monodromy computation, which is algebraic, we only know how to compute the global monodromy representation by transcendental techniques (unless in special cases). So it is probably hopeless to look for an algebraic way to compute the global monodromy exactly. On the other hand, it still makes sense to formulate the following problem:
	
	\begin{problem}
		Given $n$, $k$, and a collection of local monodromies, list all algebraic relations between the entries of the matrices comprising the monodromy representation, satisfied by all monodromy representations with prescribed local monodromies.
	\end{problem}
	
	This is the basic idea about the character variety. Before we give a precise definition we consider an example which can be worked out explicitly.
	
	\section{Fricke-Klein example}
	
	\subsection{Equations}
	This example was essentially worked out in \cite{fricke1965vorlesungen}.
	Let $k=4$, $n=2$. Let $a_1,\ldots,a_4\in\BC\setminus\{-2,2\}$. Let us try to parameterize all $4$-tuples of $2\times 2$ matrices $M_1, M_2, M_3, M_4\in \mathrm{GL}_2(\mathbb{C})$, up to conjugation, satisfying 
	\[
	\prod_{i=1}^4 M_i = I_2, \quad \det M_i = 1, \quad \trace M_i=a_i.
	\]
	Let $x=\trace (M_1 M_2)$, $y=\trace (M_1 M_3)$, $z=\trace (M_2 M_3)$. The following code computes the ideal of relations between $x,y,z$:
	
	\begin{verbatim}
		R.<m111,m112,m121,a1,m211,m212,m221,a2,m311,m312,m321,a3,x,y,z,a4>=QQ[]
		M1=matrix(2,2,[m111,m112,m121,a1-m111])
		M2=matrix(2,2,[m211,m212,m221,a2-m211])
		M3=matrix(2,2,[m311,m312,m321,a3-m311])
		eqns=[M1.det()-1,M2.det()-1,M3.det()-1,(M1*M2*M3).trace()-a4,(M1*M2).trace()-x, \
		(M1*M3).trace()-y,(M2*M3).trace()-z]
		J=R.ideal(eqns)
		J.elimination_ideal([m111,m112,m121,m211,m212,m221,m311,m312,m321])
	\end{verbatim}
	
	The output is 
	{\small
		\begin{verbatim}
			Ideal (a1*a2*a3*a4 - a1*a2*x - a1*a3*y - a2*a3*z + x*y*z - a3*x*a4 - a2*y*a4 - a1*z*a4
			+ a1^2 + a2^2 + a3^2 + x^2 + y^2 + z^2 + a4^2 - 4) of Multivariate Polynomial Ring in 
			m111, m112, m121, a1, m211, m212, m221, a2, m311, m312, m321, a3, x, y, z, a4 over 
			Rational Field	
		\end{verbatim}
	}
	
	So for each $4$-tuple of matrices satisfying the conditions, the corresponding $x, y, z\in\BC$ as above produce a solution to the so called \emph{Markov equation} with parameters $a_1, a_2, a_3, a_4$:
	\begin{multline}\label{eq:markov}
		xyz + x^2 + y^2 + z^2 - (a_1 a_2 + a_3 a_4) x - (a_1 a_3 + a_2 a_4) y - (a_1 a_4 + a_2 a_3) z \\
		+ a_1 a_2 a_3 a_4 + a_1^2 + a_2^2 + a_3^2 + a_4^2 - 4=0.
	\end{multline}
	
	Conversely, suppose $x, y, z\in\BC$ is a solution to the Markov equation. Let us see that for generically chosen values $a_1, a_2, a_3, a_4$ we can always uniquely recover the tuple of matrices, up to conjugation.
	
	First, we observe that knowing the traces $x, y, z, a_1, a_2, a_3, a_4$ we can compute the trace of an arbitrary product of matrices $M_1,\ldots,M_4$ and their inverses. Indeed, given such a product, we first replace any occurrence of $M_4$ by $M_3^{-1} M_2^{-1} M_1^{-1}$. Then using 	
	\[
	A^{-1} = \trace(A) I_2 - A \qquad (A\in \SL_2(\BC))
	\]
	we can reduce to the situation when no inverses appear. So we have an arbitrary product of $M_1, M_2, M_3$. Using
	\[
	BA = \trace(BA) I_2 - A^{-1} B^{-1} \qquad (A,B\in \SL_2(\BC))
	\]
	and getting rid of the inverses we rearrange terms modulo products of shorter length until we get to a linear combination of $I_2$, $M_1$, $M_2$, $M_3$, $M_1 M_2$, $M_1 M_3$, $M_2 M_3$, $M_1 M_2 M_3=M_4^{-1}$, all of whose traces are explicitly computed in terms of $x, y, z$ and the parameters.
	
	Knowing the traces of arbitrary products lets us find out if any $4$ products form a basis of the space of $2\times 2$ matrices. Indeed, because the trace pairing is non-degenerate, matrices $A_1,\ldots,A_{n^2}$ form a basis of the space of $n\times n$ matrices if and only if the $n^2 \times n^2$ matrix of traces $\trace(A_i A_j)$ is non-degenerate.
	
	Suppose products of $M_1,\ldots,M_4$ do span the space of matrices, i.e. the monodromy representation is irreducible, and we choose a basis $A_1, \ldots, A_4$. Knowing the traces allows us to express any other product as a linear combination of $A_1,\ldots, A_4$. In particular we recover the ring structure on the $4$-dimensional space of linear combinations of $A_1,\ldots,A_4$ and any monodromy representation provides an isomorphism of this ring with the matrix ring. Any other monodromy representation with the same traces would produce an automorphism of the matrix ring. All automorphisms of the matrix ring are inner, i.e. related by a change of basis transformation, and so we conclude that the monodromy representations with the same traces are equivalent.
	
	This was under the assumption that the monodromy representation is irreducible. If it is not, then in some basis the matrices $M_1,\ldots,M_4$ are upper triangular. In this situation we cannot determine the matrices only from the trace information. Denoting the top-left diagonal entry of $M_i$ by $\lambda_i$ we must have $\prod_{i=1}^4 \lambda_i=1$. On the other hand, $\lambda_i$ is a root of the polynomial $t^2 - a_i t + 1$ which has $2$ roots. So this kind of situation is possible only for special values of the parameters $a_1, a_2, a_3, a_4$.
	
	More precisely, we call $a_1, a_2, a_3, a_4$ \emph{generic} if none of the 16 products of the roots of $t^2 - a_i t + 1$ (we pick one root from each quadratic equation) equals $1$. We have seen, that if the parameters are generic then all monodromy representations are irreducible and in one-to-one correspondence with the solutions of the Markov equation.
	
	\subsection{Geometry of the Markov cubic}
	In algebraic geometry, given an algebraic variety $X$, i.e. a set of solutions to a system of polynomial equations, we are usually interested in its cohomology ring $H^*(X,\BC)$. Hodge theory extended by Deligne to non-compact and singular spaces \cite{deligne1971theorie, deligne1974theorie} endows the cohomology with \emph{mixed Hodge structure}. Structure constants of the ring typically encode interesting enumerative information. Sometimes the cohomology ring as a ring itself can be an interesting object of study, perhaps making a connection with a ring constructed in some other way in a completely different field of mathematics.
	
	Although the ring structure in the Fricke-Klein situation turns out to be trivial, the cohomology already has a non-trivial mixed Hodge structure. We explain how to compute it.
	
	\begin{figure}
		\begin{tikzpicture}[scale=1.]
			\draw (-1,0) -- (5,0); 
			\draw (-0.5,-1) -- (2.5,5); 
			\draw (4+0.5,-1) -- (1.5,5); 
			\draw (1.33, 0) node[draw,circle, inner sep=2pt,fill] {};
			\draw (2.67, 0) node[draw,circle, inner sep=2pt,fill] {};
			\draw (0.67, 1.33) node[draw,circle, inner sep=2pt,fill] {};
			\draw (1.33, 2.67) node[draw,circle, inner sep=2pt,fill] {};
			\draw (3.33, 1.33) node[draw,circle, inner sep=2pt,fill] {};
			\draw (2.67, 2.67) node[draw,circle, inner sep=2pt,fill] {};
		\end{tikzpicture}
		\caption{\label{fig:triangle}Markov's cubic as the complement of a triangle in the projective plane blown up at six points}
	\end{figure}
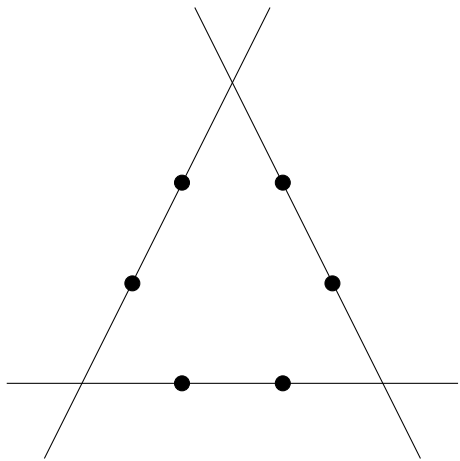
	
	Assuming the parameters $a_1,a_2,a_3,a_4\in \BC\setminus\{-2,2\}$ are generic, the cubic $X$ given by the Markov equation \eqref{eq:markov} is nonsingular. Let $\overline{X}\subset \BC\BP^3$ be the solution set of the homogenization of \eqref{eq:markov}. Its affine part is $X\subset \overline{X}$, and the complement $\overline{X}\setminus X$ is given in $\BC\BP^2$ as the solution set of the homogeneous equation $xyz=0$, thus it is a union of three lines forming a triangle. Denote this triangle by $\Delta$. We refer to Figure \ref{fig:triangle}. It is straightforward to check that homogenization of any cubic of the form 
	\[
	xyz + x^2 + y^2 + z^2 + \text{lower degree terms}
	\]
	does not have any singularities at infinity. So $\overline{X}$ is a smooth cubic surface, and therefore can be obtained as a blowup of $\BC\BP^2$ at $6$ points (for background see \cite{hartshorne1977algebraic}). The cubic $\overline{X}$, as any other smooth cubic surface contains $27$ projective lines. Three of them form $\Delta$. The intersection graph of the $27$ lines is well-known and one can see that for any triangle one can choose $6$ lines which do not intersect pairwise and do not contain any of the lines of the triangle. Each side of the triangle intersects $2$ of the $6$ lines. Blowing down these $6$ lines we obtain a map $\overline{X}\to \BC\BP^2$. The $3$ lines of $\Delta$  are isomorphically mapped to $3$ lines on $\BC\BP^2$. Choosing coordinates appropriately, one can assume that $f(\Delta)$ is the complement of $\BC^*\times\BC^*\subset \BC\BP^2$. So we have 
	\[
	U = f^{-1}(\BC^*\times\BC^*) = \BC^*\times\BC^* \subset X,
	\]
	an open subset, and the complement $X\setminus U$ is a disjoint union of $6$ copies of $\BC$:
	\[
	X = \BC^*\times\BC^* \sqcup \underbrace{\BC \sqcup \cdots \sqcup \BC}_{\text{$6$ copies}}.
	\]  
	The $6$ copies of $\BC$ are divided into $3$ pairs according to which of the three lines of $\Delta$ they intersect.
	
	\subsection{Cohomology}
	With this picture in mind, we have several ways to compute the cohomology. All cohomology groups we consider have complex coefficients. From the relative cohomology sequence for the pair $\BC^*\times\BC^* = U\subset X$ we obtain an exact sequence
	\[
	0\to H^1(X) \to H^1(U) \to H^2(X,U) \to H^2(X) \to H^2(U) \to 0.
	\] 
	Since both $\BC$ and $X$ are non-singular, by excision we have that the only non-vanishing relative cohomology group is $H^2(X, U)\cong \BC^6$. So the remaining part of the long exact sequence vanishes, except for the isomorphism $H^0(X)\cong H^0(U)\cong \BC$. The map $\BC^2 = H^1(U) \to H^2(X,U)=\BC^6$ is dual to the map $H_2(X,U) \to H_1(U)$. Here a point class on $\BC$ corresponds to a disc class in $H_2(X,U)$ whose boundary is a $1$-cycle in $U$ which gets contractible in $U\cup \BC$. From the picture we see that all $1$-cycles become contractible, so the map $H_2(X,U) \to H_1(U)$ is surjective, the dual map is injective and we have $H^1(X)=0$, $H^2(X)$ is an extension of $H^2(U)=\BC$ by $H^2(X,U)/H^1(U)\cong\BC^4$.
	
	This is reflected by the weight filtration on $H^2(X)$.  The cohomology $H^2(X,U)/H^1(U)\subset H^2(X)$ as a mixed Hodge structure is pure of weight $2$, while $H^2(U)$ has weight $4$. Thus we have $W_2 = H^2(X,U)/H^1(U)\subset H^2(X)$ with $\dim W_2=4$, $W_3=W_2$, and $W_4 = H^2(X)$ with $\dim W_4/W_2 = 1$.
	
	Another way to compute the cohomology is to use the compactification. Let us denote $X_0=\overline{X}$, $X_1 = X_0\setminus \BC\BP^1$, $X_2=X_0\setminus (\BC\BP^1 \cup \BC\BP^1)$, $X_3=X_0\setminus \Delta$, so that we pass from $X_i$ to $X_{i+1}$ by removing a line. This leads to exact sequences for each $i$:
	\[
	0\to H^1(X_{i+1})\to H^2(X_i, X_{i+1}) \to H^2(X_i)\to H^2(X_{i+1}) \to H^3(X_i,X_{i+1}) \to 0
	\]
	\[
	0\to H^3(X_{i+1}) \to H^4(X_i,X_{i+1}) \to H^4(X_i)\to H^4(X_{i+1}) \to 0,
	\]
	as we will see that the even cohomologies of $X_i$ vanish at each step. At each step we will use that $H^k(X_i, X_{i+1})=H^{k-2}(X_i\setminus X_{i+1})$ by excision, which is non-zero only for $k=0,1,2$.
	Initially we have $H^0(X_0)=\BC$, $H^2(X_0)=\BC^7$ (represented by $7$ algebraic cycles: the line in $\BC\BP^2$ and the $6$ exceptional curves), $H^4(X_0)=\BC$, all other cohomologies vanish. 
	
	For $i=0$ we have $X_0\setminus X_1=\BC\BP^1$, which has $H^0=\BC$ and $H^2=\BC$. These coincide with $H^2(X_0,X_1)$ and $H^4(X_0,X_1)$ respectively and we see that in $H^2(X_1)$ and in $H^4(X_1)$ we have one less dimension than in $H^2(X_0)$, respectively $H^4(X_0)$. So $H^4(X_1)=0$, $H^2(X_1)=\BC^6$.
	
	For $i=1$ we have $X_1\setminus X_2 = \BC$, which only has $H^0=\BC$. So we have that $X_2$ has one dimension less in $H^2(X_2)$, nothing else changes: $H^2(X_1)=\BC^5$.
	
	Finally, for $i=2$ we have $X_2\setminus X_3 = \BC^*$, which has $H^0=\BC$ and $H^1=\BC$, contributing to $H^2(X_2,X_3)$ and $H^3(X_2, X_3)$ correspondingly. Thus two things happen. First, only $4$ dimensions from $H^2(X_2)$ out of $5$ contribute to $H^2(X_3)$. On the other hand, one extra dimension comes from $H^3(X_3,X_2)$. The net total is $H^2(X_3)=\BC^5$, as in $X_2$, but the nature of the $5$ dimensions is different, and correspondingly the weights are different. The $4$ dimensions coming from $H^2(X_2)$ have weight $2$, and this is reflected by the fact that they are represented by algebraic cycles. The remaining class has weight $4$ as it comes from $H^3(X_3,X_2)$. The weight of $H^1(\BC^*)$ is $2$, but it gets increased by $2$ when we pass to $H^3(X_3,X_2)$. 
	
	These findings can be summarized in the table on Figure \ref{fig:cohomology}.
	
	\begin{figure}
		\begin{tabular}{c|c|c}
			$H^*$ & $W_i H^i$ & $W_{i+2} H^i / W_{i} H^i$ \\
			\hline
			$H^2$ & $\mathbb{C}^4 = \langle c_1(L_i)\;|\;i=1,2,3,4\rangle$ & $\mathbb{C} = \langle \omega=\text{symplectic form} \rangle$ \\
			$H^0$ & $\mathbb{C} = \langle 1 \rangle$ &
		\end{tabular}
		\caption{\label{fig:cohomology}Cohomology of the Markov cubic}
	\end{figure}
	
	The $4$ algebraic classes in $H^2(X)$ can be represented by characteristic classes of line bundles. The extra class can be represented by an explicit holomorphic symplectic form as explained below. Using the natural coordinates on $\BC^*\times \BC^*$, the form is written as
	\[
	\omega = \frac{dx}x \wedge \frac{dy} y,
	\]	
	and one can verify that the form extends to the whole of $X$.
	
	\subsection{Mixed Hodge Polynomial}\label{ssec:mixed hodge polynomial}
	From Figure \ref{fig:cohomology} we obtain that the Poincar\'e polynomial of $X$ is $1 + 5 q^2$ and the mixed Hodge polynomial is $P(q,t) = 1 + 4 q^2 t^2 + q^2 t^4$. We prefer to consider a slightly different version of the mixed Hodge polynomial, without losing any information. Let
	
	\[
	\mathrm{MH}(X; q,t) := q^{\frac12 \dim_\BC X} P(X;q^{-\frac12} t^{-\frac12}, t^{\frac12}).
	\]
	For instance, in our running example we obtain
	\[
	\mathrm{MH}(X; q,t) = q + t + 4.
	\]
	
	We should also compare our calculations with the number of points on our variety over $\BF_q$. Assuming all the lines are defined over the base field, we see that the first computation produces
	\[
	(q-1)^2 + 6 q  = q^2 + 4 q + 1.
	\]
	The second computation produces
	\[
	q^2 + 7q + 1 - (q+1) - q - (q-1) = q^2 + 4 q + 1.
	\]
	In each case, one can track how the different powers of $q$ contribute to the cohomology so that $q^2$ corresponds to $H^0$ and $4q + 1 $ corresponds to $H^2$.
	
	In any case, for a general algebraic variety $X$, if the number of points is a polynomial of $q$, we have
	\begin{equation}\label{eq:Fq points and mixed Hodge}
		|X(\mathbb{F}_q)| = q^{\dim_\BC X} P(X;-1,q^{-\frac12}) = q^{\frac12 \dim_\BC X} \mathrm{MH}(X;q,q^{-1}),
	\end{equation}
	assuming $\mathrm{MH}$ is indeed a polynomial in $q$, $t$, which means that all the weights are even and the cohomology is concentrated in even degrees.
	
	\section{Setup}
	After reviewing in detail the basic example, we turn to the general setup and precise definitions. There is no reason to consider only genus $0$ surfaces, and the classical motivation directly extends to the more general case, for instance ODEs can be studied on complex curves of arbitrary genus.
	
	\begin{definition}
		Fix $g\geq 0$, $n, k\geq 1$ and Jordan forms $J_1, J_2, \ldots, J_k \in\mathrm{GL}_n$. The associated genus $g$ \emph{character variety} $X=\mathrm{Char}_{g;J_1,\dots,J_k}$ is the space
		\[
		\{(A_1,\ldots,A_g,B_1\ldots,B_g,M_1, \cdots, M_k): M_i\sim J_i,
		\]
		\[ A_1 B_1 A_1^{-1} B_1^{-1} \cdots A_g B_g A_g^{-1} B_g^{-1} M_1 \cdots M_k = I\}/\sim,
		\]
		where $\sim$ identifies tuples of matrices which are conjugate with respect to the $\GL_n(\BC)$ action.
	\end{definition}

	\begin{definition}[\cite{hausel2011arithmetic}]\label{def:generic}
		The data $J_1,\ldots J_k$ are \emph{generic} if 
		\begin{enumerate}
			\item $\prod_{i=1}^k \det J_i=1$
			\item $\prod_{i=1}^k \det J_i'\neq 1$ for any $0<n'<n$ and any block-decomposition $\left(\begin{array}{c|c} J_i' & J_i''\\ \hline 0 & J_i'''\end{array}\right)\sim J_i$ ($J_i'\in \mathrm{GL}_{n'}$).
		\end{enumerate}
	\end{definition}
	
	\begin{theorem}[\cite{hausel2011arithmetic}]
		For generic data $n, k, J_1, \ldots, J_k$, the quotient space $X=\mathrm{Char}_{g;J_1,\dots,J_k}$ has the structure of a smooth affine algebraic variety.
	\end{theorem}
	In fact, considering functions on $X$ of the form
	\[
	\trace(M_{i_1}\text{or} A_{i_1} \text{or} B_{i_1} \cdots M_{i_m}\text{or} A_{i_m} \text{or} B_{i_m})
	\]
	we can realize $X$ as an algebraic subset of $\mathbb{C}^N$ for suffuciently large $N$.
	
	Basically, the same reasoning as in the Fricke-Klein situation can be applied in general. However, the number of traces one needs to consider is uncontrollably large, and it seems hopeless to be able to produce all equations by the elimination computation. Even if we do obtain a complete list of equations, it is not clear which interesting information we should be able to extract from such a list. For instance, equations in the case of genus $0$, rank $2$, $5$ punctures are given in \cite{simpson2017explicit}. We do hope, however, that computing the cohomology ring is an achievable goal.
	\begin{problem}
		Compute the cohomology ring $H^*(X,\mathbb{C})$.
	\end{problem}
	
	\section{What is known?}
	Having defined our basic object of study, we explain some important structural results.
	
	\subsection{Non-abelian Hodge correspondence}
	The non-abelian Hodge correspondence, pioneered by Hitchin \cite{hitchin1987self} and developed by Simpson \cite{simpson1990harmonic, simpson1991nonabelian, simpson1992higgs}, is a construction of a remarkable bijection between local systems and stable Higgs bundles, which induces a diffeomorphism between the character variety and the moduli space that parametrizes stable Higgs bundles.
	
	To define Higgs bundles, we need to view the Riemann surface with punctures as a compact algebraic curve $C$ over $\BC$ with marked points. Then a \emph{Higgs bundle} is a pair $(\mathcal{E}, \theta)$, where $\mathcal{E}$ is an algebraic vector bundle on $C$ and $\theta$ a \emph{Higgs field}, which is a twisted matrix-valued function on $C$ of a certain kind. More precisely, $\theta:\mathcal{E}\to\mathcal{E}\otimes\omega_C(D)$, where $D=\sum_{i=1}^k [t_i]$ is the divisor on $C$ corresponding to the marked points. Explicitly, this means that on an open subset of $C$ with a trivialization of $\mathcal{E}$ is trivial, $\theta$ is a matrix whose entries are algebraic differential forms on $C$ with simple residues at the marked points.
	
\begin{figure}
		\begin{tabular}{c|c}
			Local systems  & Higgs bundles  \\
			\hline
			$\pi_1(C)\to \mathrm{GL}_n(\mathbb{C})$ & $(\mathcal{E}, \theta)$: $\mathcal{E}$ algebraic vector bundle on $C$, \\
			& $\theta:\mathcal{E}\to\mathcal{E}\otimes\omega_C$ Higgs field \\
			
			punctures, $J_i$ & parabolic structure, stability \\
			
			$X_B$ ``Betti'' moduli space & $\cong X_D$ ``Dolbeault'' moduli space \\
			
			\begin{tikzpicture}[scale=0.5]
				\draw (-1,0) -- (5,0); 
				\draw (-0.5,-1) -- (2.5,5); 
				\draw (4+0.5,-1) -- (1.5,5); 
				\draw (1.33, 0) node[draw,circle, inner sep=2pt,fill] {};
				\draw (2.67, 0) node[draw,circle, inner sep=2pt,fill] {};
				\draw (0.67, 1.33) node[draw,circle, inner sep=2pt,fill] {};
				\draw (1.33, 2.67) node[draw,circle, inner sep=2pt,fill] {};
				\draw (3.33, 1.33) node[draw,circle, inner sep=2pt,fill] {};
				\draw (2.67, 2.67) node[draw,circle, inner sep=2pt,fill] {};
			\end{tikzpicture}
			&
			\includegraphics[width=4.5cm]{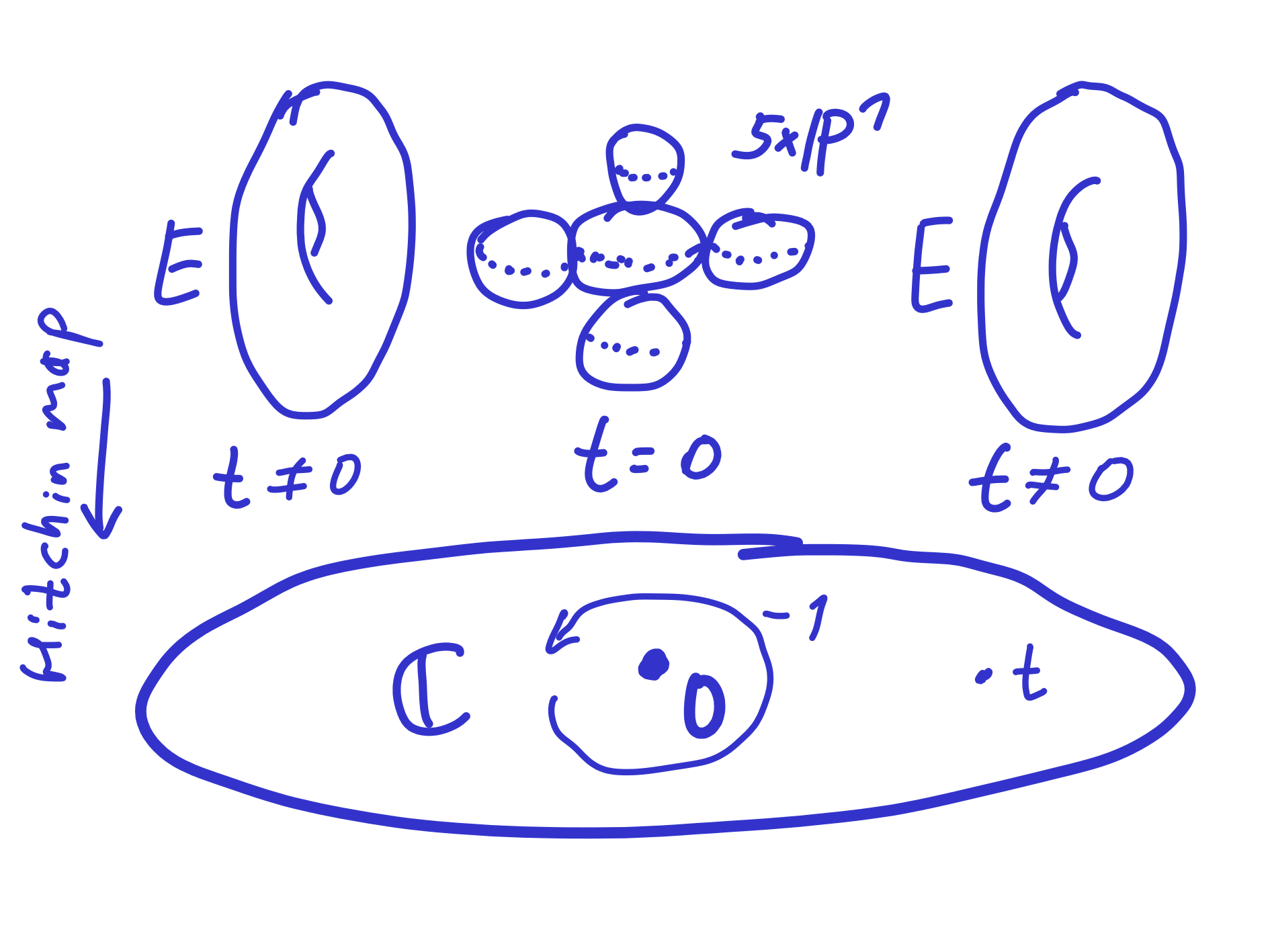} \\
			
			affine algebraic variety & Hitchin system $X_D\to \mathbb{C}^{\text{dimension}/2}$\\
			
			Cohomology ring $H^*(X_B,\mathbb{C})$ & the same ring \\
			
			Weight ``W'' filtration $W_* H^*$ & perverse ``P'' filtration\\
		\end{tabular}
		\caption{\label{fig:nonabelian hodge}Non-abelian Hodge correspondence}
	\end{figure}

	Furthermore, these Higgs bundles are required to be endowed with \emph{parabolic structures} at the marked points. This is simply a choice of a partial flag at the fiber of $\mathcal{E}$ at each marked point, preserved by the residue $\res_{t_i} \theta$, which acts on that fiber by a linear transformation. The steps of the partial flag precisely match the dimensions of the generalized eigenspaces of the Jordan forms $J_i$. The data of a Higgs bundle with a parabolic structure has to satisfy a \emph{stability condition}. The data of the eigenvalues and the Jordan blocks of $J_i$ translate into \emph{parabolic weights} defining the stability condition, and into the Jordan form of $\res_{t_i} \theta$ acting on the subquotients of the flag. Details are given in \cite{simpson1990harmonic}. We summarize the features appearing on both sides of this correspondence in the table on Figure \ref{fig:nonabelian hodge}.
	
	\subsection{Generators of $H^*(X)$}
	Generators of the cohomology ring are known. There is a general argument in \cite{markman2002generators} concerning moduli spaces of sheaves on surfaces. With a little modification it can be applied to the situation of parabolic Higgs bundles. Basically speaking, certain Chern class of the complex of vector bundles on $X_D\times X_D$ computing the derived $\Hom$ between Higgs bundles is a multiple of the diagonal class. Moreover, there exists a compactification $\overline{X_D}\supset X_D$ such that the complex extends to $X_D\times \overline{X_D}$ and the Chern class remains a multiple of the diagonal class. Thus the identity map $H^*(X_D)\to H^*(X_D)$ can be written in terms of the tautological classes.
	
	Consider the tautological bundle\footnote{In fact the tautological bundle often is not well-defined as a vector bundle. It is well-defined as a projective bundle, and the ``determinant'' $\det \CE$ is also well-defined, but the determinant might have no $n$-th root and correspondingly $c_1(\det \CE)\in H^2(X_D,\BZ)$ may not be divisible by $n$. Since we are dealing with the cohomology with complex coefficients, this does not create a problem and all the Chern classes are still well-defined.} $\mathcal{E}\to X_D\times C$. For each $i\in\{0,1,2\}$, $\xi\in H^i(C)$ and $m\in\BZ_{\geq 0}$ define
	\begin{equation}\label{eq:definition e}
		e_m(\xi):=(\pi_X)_*( c_m(\mathcal{E}) \cup \xi)\in H^{2m+i-2}(X_D).
	\end{equation}
	Here $c_m(\mathcal{E})\in H^{2m}(X_D\times C)$ is the Chern class and $\pi_{X}$ is the projection $X_D\times C\to X_D$. If we have a marked point $t_j$, then each point of $X_D$ gives a flag $0=F_{j0}\subset F_{j1} \subset\cdots\subset F_{j r_j} = \CE(t_j)$. Each subquotient of the flag produces a vector bundle on $X_D$ and we take its Chern classes $c_{m,j,\ell} := c_m(F_{j \ell} / F_{j \ell-1})$.
	\begin{theorem}
		The cohomology ring $H^*(X_D)$ is generated by the classes $e_m(\xi)$ and $c_{m,j,\ell}$.
	\end{theorem}
	
	The cohomology ring of the Riemann surface $C$ has a basis $1\in H^0(C)$, $\mathrm{pt}\in H^2(C)$ and $\sigma_1,\ldots,\sigma_{2g}\in H^1(C)$ with the product structure given by $\sigma_i \sigma_{i+g} = -\sigma_{i+g}\sigma_i=\mathrm{pt}$, all other $\sigma_i \sigma_j$ being zero.
	
	\subsection{Explicit relations of $H^*(X)$} are known for $n=2$, $k=1$, $J_1=-I$, arbitrary $g$ (\cite{hausel2003relations}). The cohomology is generated by certain elements $\varepsilon_1,\ldots,\varepsilon_{2g}\in H^1(X_D)$, $\alpha\in H^2(X_D)$, $\psi_1,\ldots,\psi_{2g}\in H^2(X_D)$ and $\beta\in H^4(X_D)$. These are certain multiples of $e_1(\sigma_i)$, $e_2(1)$, $e_2(\sigma_i)$ and $e_2(\mathrm{pt})$. The elements $\varepsilon_i$ generate a copy of the cohomology ring of the Jacobian of $C$ and we have
	\[
	H^*(X_D) = H^*(\Jac(C)) \otimes \bigoplus_{i=0}^g \wedge_0^i(\psi) \otimes \BC[\alpha,\beta,\gamma] / J_{i}^{g-i},
	\]
	where $J_{i}^{g-i}\subset \BC[\alpha,\beta,\gamma]$ is the ideal generated by an explicit list of polynomials in $\alpha, \beta, \gamma$. Here $\gamma=-2 \sum_{i=1}^g \psi_i \psi_{i+g}$ and $\wedge_0^i(\psi)$ denotes the primitive part in the space of $i$-fold products of $\psi_j$-classes. 
	
	\subsection{Weight filtration} The weight filtration is very easily described in terms of the tautological generators. The following was shown in \cite{shende2016weights} in the case without parabolic points and in \cite{mellit2019cell} in a more general case:
	\begin{theorem}
		Define the weight of the tautological generators $e_{m}(\xi)$ and $c_{m,j,\ell}$ to be $2 m$. Then $W_m H^*(X_B)$ is spanned by products of tautological generators of total weight $\leq m$.
	\end{theorem}
	
	\subsection{The pure part} 
\begin{figure}
		\includegraphics[width=.5\linewidth]{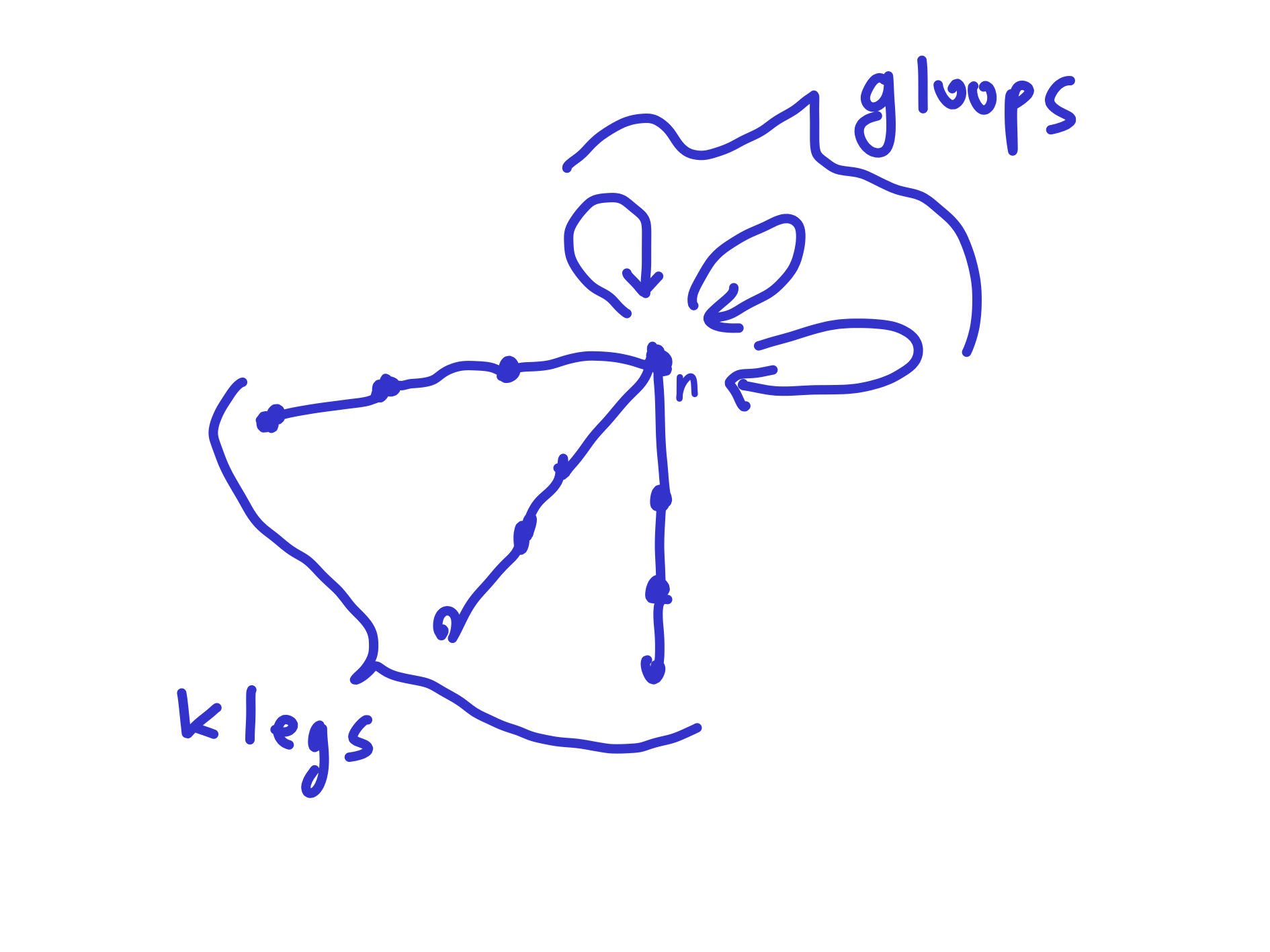}
		\caption{\label{fig:comet} Comet-shaped quiver with $k$ legs and $g$ loops.}
	\end{figure}

The \emph{pure part} $H^*_{\mathrm{pure}}(X_B)=\bigoplus_i W_i H^i(X_B)$ is naturally a subring of $H^*(X_B)$. In \cite{hausel2011arithmetic} it was conjectured that the Poincar\'e polynomial of $H^*_{\mathrm{pure}}(X_B)$ matches that of the quiver variety associated to the corresponding comet-shaped quiver, see Figure \ref{fig:comet}.
	
	\subsection{$\boldsymbol{P}=\boldsymbol{W}$} Computing the characteristic polynomial of $\theta:\CE\to \CE\otimes \omega_C(D)$ at all points of $C$ produces a well-defined polynomial map $X_D\to \BC^{\dim_\BC X/2}$ called the \emph{Hitchin map}, which is known to be proper. This induces a filtration $P_* H^*(X_D)$ on $H^*(X_D)$ called the \emph{perverse Leray filtration}  (\cite{cataldo2009decomposition}). According to the recently proved $P=W$ conjecture (\cite{cataldo2012topology}), we have
	\begin{theorem}[\cite{maulik2024pw}, \cite{hausel2022p}]
		\[
		P_i=W_{2i}.
		\]
	\end{theorem}
	
	\subsection{$\boldsymbol{q}-\boldsymbol{t}$ symmetry} Recall the Mixed Hodge polynomial $q+t+4$ from Section \ref{ssec:mixed hodge polynomial}. Its symmetry in the variables $q$ and $t$ was called the \emph{curious hard Lefschetz} property in \cite{hausel2011arithmetic}, and it was expected to hold in general. Having $P=W$, it follows from the \emph{relative hard Lefschetz theorem} on the $X_D$-side. We have
	
	\begin{theorem}[\cite{mellit2019cell}]
		$H^*(X_B)$ satisfies curious hard Lefschetz:
		\[
		\omega^i: \mathrm{Gr}^W_{\dim_\BC X-2i} H^j(X_B) \to \mathrm{Gr}^W_{\dim_\BC X+2i} H^{j+2 i}(X_B)
		\]
		is an isomorphism.
	\end{theorem}
	
	In fact, the proof of $P=W$ from \cite{maulik2024pw} relies on this result. On the other hand, the proof of $P=W$ in \cite{hausel2022p} is independent of it, and therefore provides another proof of the curious hard Lefschetz.
	
	Similar curious hard Lefschetz property holds for a certain version of the Khovanov homology of knots (conjectured in \cite{dunfield2006superpolynomial}, proved in \cite{gorsky2024tautological}).
	
	\subsection{Partition function} Perhaps the biggest problem in the field that still remains open is the conjectural formula for the generating function of the mixed Hodge polynomials of all character varieties given in \cite{hausel2011arithmetic} in terms of a partition function. Here $k$ (the number of marked points) and $g$ are fixed. For each $i=1,\ldots,k$ let $\mathbf{x}_i:=(x_{i1},x_{i2},\ldots)$ be an infinite sequence of variables. Suppose we have a sequence $\boldsymbol\mu_i:=(\mu_{i1},\mu_{i2},\ldots)$ of non-negative integers with $\sum_{j}\mu_{ij}=n$. To this we associate a conjugacy class $J_{\mu_i}$ of the diagonal matrix with eigenvalues having multiplicities $\mu_{i1},\mu_{i2},\ldots$. For a tuple of such sequences $\boldsymbol{\mu}=(\mu_1,\ldots,\mu_k)$ we choose the eigenvalues arbitrarily so that the conjugacy classes satisfy the genericity condition (Definition \ref{def:generic}) and consider the corresponding character variety $X_{g;J_{\mu_1},\cdots,J_{\mu_{k}}}$. We take its mixed Hodge polynomial (Section \ref{ssec:mixed hodge polynomial}) with weight
	\[
	\mathbf{x}^{\boldsymbol\mu}:=\prod_{i=1}^k \prod_j x_{ij}^{\mu_{ij}}.
	\]
	The result is a symmetric function in each sequence of variables $\mathbf{x}_i$ whose coefficients are polynomials in $q$ and $t$:
	\[
	\mathrm{MH}_{g,k}(\mathbf{x}_1,\ldots,\mathbf{x}_k; q,t) :=
	\sum_{n=1}^\infty \sum_{ \mu_1, \dots, \mu_{k} \vdash n} \mathrm{MH}(X_{g;J_{\mu_1},\cdots,J_{\mu_{k}}};q,t) \; \mathbf{x}^{\boldsymbol\mu}.
	\]
	On the other hand, we define the \emph{partition function}
	\begin{equation}\label{eq:partition function}
		\Omega_{g,k}(\mathbf{x}_1,\ldots,\mathbf{x}_{k}; q,t) 
		:=\sum_{n=0}^\infty \sum_{\lambda\vdash n} \tilde H_\lambda(\mathbf{x}_1;q,t) \cdots \tilde H_\lambda(\mathbf{x}_k;q,t)
		\prod_{\square\in\lambda}\frac{\left(q^{a(\square)+1/2} - t^{\ell(\square)+1/2}\right)^{2g}}{(q^{a(\square)} - t^{\ell(\square)+1})(q^{a(\square)+1} - t^{\ell(\square)})},
	\end{equation}
	where $\lambda$ goes over all partitions (Young diagrams), $\tilde H_\lambda(\mathbf{x}_k;q,t)$ denotes the \emph{modified Macdonald polynomial}, the product goes over the cells of the Young diagram, and $a,\ell$ denote the usual arms and legs statistics (see \cite{garsia1996remarkable}).
	
	\begin{conjecture}[\cite{hausel2011arithmetic}]\label{conj:partition function conjecture}
		\begin{equation}\label{eq:partition function conjecture}
			\Omega_{g,k} = \pExp\;\left(-\frac{\mathrm{MH_{g,k}}}{(q-1)(t-1)}\right).
		\end{equation}
	\end{conjecture}
	
	Here $\pExp$ denotes the \emph{plethystic exponential}. It can be defined for any expression involving some variables $a,b,c,\ldots$ as follows:
	\[
	\pExp(f(a,b,c,\ldots)) := \exp\left(\sum_{n=1}^\infty \frac{f(a^n,b^n,c^n,\ldots)}{n}\right).
	\]
	If $a,b,c,\ldots$ are monomials in some variables, then we have
	\[
	\pExp(a+b+c+\cdots) = \frac{1}{1-a} \; \frac{1}{1-b} \; \frac{1}{1-c} \; \cdots,
	\]
	so plethystic exponential can be thought of as the conversion between infinite sums and infinite products. There is an inverse operation $\pLog$ called the \emph{plethystic logarithm}. So if we focus on computing the mixed Hodge polynomials from \eqref{eq:partition function conjecture}, we can write
	\begin{equation}\label{eq:partition function logarithm}
		\mathrm{MH_{g,k}} = -(1-q)(1-t) \pLog \Omega_{g,k}.
	\end{equation}
	
	\subsection{What is known?}
	Conjecture \eqref{conj:partition function conjecture} was motivated by a computation performed in \cite{hausel2011arithmetic}. It was shown that its implication for the number of $\mathbb{F}_q$ points is correct (see \eqref{eq:Fq points and mixed Hodge}), by counting representations $\pi_1(C)\to \mathrm{GL_n}(\mathbb{F}_q)$ using the character theory of the finite group $\mathrm{GL_n}(\mathbb{F}_q)$.
	\begin{theorem}[\cite{hausel2011arithmetic}]
		The number of points over $\mathbb{F}_q$ ($t=q^{-1}$) deduced from \eqref{eq:partition function conjecture} is correct.
	\end{theorem}
	
	Initially, in the identity \eqref{eq:partition function logarithm} it was not at all clear that the right hand side has coefficients which are polynomials in $q$, $t$, which they should be because the mixed Hodge polynomials on the left hand side are polynomials by definition.	
	Using a new identity in Macdonald polynomials combinatorics established in \cite{garsia2019five}, the author first proved
	\begin{theorem}[\cite{mellit2018integrality}]
		The coefficients of $(q-1)(t-1)\; \mathrm{Log}\;\Omega$ are polynomials in $q$ and $t$.
	\end{theorem}
	
	Then, together with an approach to count Higgs bundles over $\mathbb{F}_q$ by Schiffmann (\cite{schiffmann2016indecomposable}, \cite{mozgovoy2020counting}) and a new relationship between the Macdonald polynomials and the affine Springer fibers, the author proved
	\begin{theorem}[\cite{mellit2020poincare}, \cite{mellit2020poincarea}]\label{thm:poincare}
		The Poincar\'e polynomials ($t=1$) deduced from \eqref{eq:partition function conjecture} are correct.
	\end{theorem}
	
	With the help of this theorem, there is a following line of attack on Conjecture \eqref{conj:partition function conjecture}. Suppose we managed to construct a ``candidate model'' $M$ for the cohomology ring of the character variety so that it comes with a surjective map $M\to H^*(X)$. Showing surjectivity of the map is plausible because the generators of $H^*(X)$ are known. If we furthermore manage to show that the dimension of $M$ is exactly as predicted by Conjecture \eqref{conj:partition function conjecture}, then by Theorem \ref{thm:poincare} the map $M\to H^*(X)$ will have to be an isomorphism. Partial progress in this direction is discussed in the next sections.

	\section{COHAs/Coulomb branches}
	
	In order to explain our approach we need to introduce some ideas.
	\subsection{Cohomological Hall Algebra (COHA)}
	Denote by $\mathscr{H}_2$ the Lie algebra of Poisson vector fields on $\mathbb{C}^2$. Explicitly, it is the Lie algebra with basis elements $T_{r,s}$ for each pair $r,s\geq 0$, corresponding to the vector field with Hamiltonian $x^r y^s$. It is convenient not to impose $T_{0,0}=0$, although the corresponding vector field vanishes. The Lie bracket is given as follows:
	\[
	[T_{r,s}, T_{r',s'}] = (r's - r s') T_{r+r'-1,s+s'-1}\qquad(r,s\geq 0).
	\]
	The tensor product $\mathscr{H}_2 \otimes H^*(C)$ is naturally a super-Lie algebra. Let $H^*(X)[x,y]$ be the cohomology ring with two formal variables $x$, $y$ adjoined. Let $n$ and $d$ be the rank and the degree of the Higgs bundles parameterized by the moduli space $X$, and let $\chi:=d + (1-g)n$ be the Euler characteristic.
	\begin{theorem}[\cite{hausel2022p}, \cite{mellit2023coherent}]
		The Lie algebra $\mathscr{H}_2 \otimes H^*(C)$ acts on $H^*(X)[x,y]$ in such a way that $T_{0,m}\otimes \xi$ satisfies 
		\[
		\left((T_{0,m}\otimes \xi) f\right)(0,0) = \psi_m(\xi)(f) := \left(p_m(\xi) + m(1-g) p_{m-1}(\xi\cup \mathrm{pt})\right)\cup f \qquad (f\in H^*(X))
		\]
		and $T_{1,0}\otimes \mathrm{pt}$, $T_{1,0}\otimes 1$, $T_{0,1}\otimes \mathrm{pt}$, $T_{0,1}\otimes 1$ act as $x$, $-n \frac{\partial}{\partial y}$, $y$, $n \frac{\partial}{\partial x} + \chi$ respectively.
	\end{theorem}
	Here $p_m(\xi)$ is defined similarly to $e_m(\xi)$ in \eqref{eq:definition e}, except that we use the power sums of the Chern roots, i.e. $m!$ times the $m$-th component of the Chern character. Thus for instance, the central operator $T_{0,0}\otimes\mathrm{pt}$ acts via the multiplication by the rank of the Higgs bundles $n$. The operator $T_{0,1}\otimes\mathrm{1}$ commutes with $y$ and acts on elements of $H^*(X)[y]$ via the multiplication by $d+(1-g)n = \chi$.

	Note that quadratic Hamiltonians correspond to the linear action of $SL(2,\BC)$ on $\BC^2$, and so the elements $T_{2,0}$, $T_{1,1}$ and $T_{0,2}$ span a copy of $\mathfrak{sl}_2\subset \mathscr{H}_2$. In fact it is this $\mathfrak{sl}_2$ that was matched with the relative Lefschetz $\mathfrak{sl}_2$ in \cite{hausel2022p} in order to pin down the perverse filtration algebraically.
	
	\subsection{Coulomb branch}
	The construction described in \cite{braverman2018towards} associates to a pair of a reductive group $G$ and its representation $N$ an affine Poisson algebraic variety called the Coulomb branch. We are interested in the situation where $G=GL_n\times \BC^*$ acts on the space of matrices $N=\BC^{n\times n}$: $GL_n$ acts by conjugation and $\BC^*$ by the scalar multiplication. Then Higgs bundles can be seen as maps from the curve $C$ to the stack $N/G$. Our hope is that the Coulomb branch has something to do with the cohomology of our moduli spaces. This is the situation of the adjoint representation with flavor and the Coulomb branch is related to the Hilbert scheme (\cite{braverman2021line}) of $\BC^*\times \BC$, which is a symplectic resolution of singularities of $(\BC^*\times \BC)^n / S_n$. It is not completely clear why we should pass to the Hilbert scheme of $\BC^2$, but our construction in \cite{hausel2022p}, in fact first produces an algebra of Hamiltonian fields on $\BC^*\times\BC$, and only then it is degenerated to the algebra for $\BC^2$. Hopefully, this passage will be clarified by experts in Coulomb branches.
	
	Now suppose we are on $\mathrm{Hilb}^n(\BC^2)$. Then $\mathscr{H}_2$ acts on many natural objects such as functions, differential forms etc. So even if we do not completely understand the Coulomb branch picture, it is a natural guess to look for our model of the cohomology ring among natural objects living on the Hilbert scheme.
	
	\subsection{Geometric engineering}
	An indication similar to the one above is given in \cite{chuang2015parabolic}, where a physical explanation of Conjecture \ref{conj:partition function conjecture} is provided via interpreting the left hand side (the partition function) as
	\begin{equation}\label{eq:geometric engineering}
		\Omega = \chi(\mathrm{Hilb}_n(\mathbb{C}^2), P^{\otimes k} \otimes (\Omega^{*})^{\otimes g}\otimes O(g-1)).
	\end{equation}
	Here $P$ is the Procesi bundle ($S_n$-equivariant bundle coming from Haiman's identification $\mathrm{Hilb}_n(\mathbb{C}^2)=S_n-\mathrm{Hilb}(\mathbb{C}^{2n})$), $\Omega^*$ is the bundle of differential forms, $\chi$ denotes the equivariant Euler characteristic.
	
	\subsection{Proposal}
	One problem with \eqref{eq:geometric engineering} is that the right hand side does not always have positive coefficients. So the vector bundle on the right hand side of \eqref{eq:geometric engineering} has higher cohomologies. Thus there seem to be no way to upgrade the identity \eqref{eq:partition function conjecture} to an isomorphism between vector spaces using \eqref{eq:geometric engineering}. One can trace the problem to the fact that the factor $-\frac{1}{(1-q)(1-t)}$ is manifestly not positive, when expanded in $q$ and $t$. To resolve this issue we propose to expand in $q^{-1}$ and $t$ instead. Geometric meaning of this expansion is different from taking the global sections. Instead, we should take the local cohomology with respect to a Lagrangian subvariety.
	
	\begin{theorem}[M.-Romero, in preparation]
		Let $L\subset \mathrm{Hilb}_n(\mathbb{C}^2)$ be the Lagrangian cut out by the equations $\sum_{i=1}^n x_i^r=0$ ($r>0$). Let $\mathcal{E}$ be any vector bundle on $\mathrm{Hilb}_n(\mathbb{C}^2)$. Then for any $i\neq n$ we have
		\[
		H^i_L(\mathrm{Hilb}_n(\mathbb{C}^2), \mathcal{E})=0.
		\]
	\end{theorem}
	Here $H^i_L$ is the functor of local cohomology (see \cite{hartshorne1967local}), the derived functor to the functor of taking sections supported at $L$. Then we expect that \eqref{eq:partition function conjecture} has the following interpretation:
	\begin{conjecture}
		$
		H^n_L(P^{\otimes k} \otimes (\Omega^{*})^{\otimes g}\otimes O(g-1))\;\check{}\; \cong H_*^{\mathrm{BM}}(\text{Higgs moduli stack}).
		$
	\end{conjecture}
	Here it is not clear which precise Higgs moduli stack should be taken on the right hand side. The whole Higgs moduli stack has infinitely many irreducible components and the expectation is that we need to take a nice open substack containing only finitely many components in each degree, possibly imposing a bound on the slope. For a bigraded vector space $M$ we denote by $M\;\check{}\;$ the bigraded dual vector space.
	
	\section{Main example} A special case can be worked out in full detail. We take $g=0$, i.e. $C=\BP^1$, and consider the $k+1$-twisted Higgs bundles, so that the Higgs field is a morphism $\CE\to \CE(k-1)$. Equivalently, we may consider $k+1$ marked points without any parabolic structures and conditions on the residue of the Higgs field. The corresponding specialization of the partition function is given by
	\begin{equation}\label{eq:omega example}
		\Omega_k(z; q,t) := 
		\sum_{n=0}^\infty \sum_{\lambda\vdash n} \frac{(-1)^{(k+1)|\lambda|}\left(q^{n(\lambda')} t^{n(\lambda)}\right)^{k+1}}{
			\prod_{\square\in\lambda} (q^{a(\square)} - t^{\ell(\square)+1})(q^{a(\square)+1} - t^{\ell(\square)})} z^{n}.
	\end{equation}
	We use the variable $z$ to keep track of the rank. In \eqref{eq:partition function} we didn't need it because the rank of each term can be extracted as the degree of the symmetric function in each set of variables.
	The first terms of $\Omega$ look as follows:
	\[
	1 + \frac{(-1)^{k}}{(1-t)(1-q)} z + \frac{q^{k+1}(1-t^2)-t^{k+1}(1-q^2)}{(q-t)(1-q)(1-q^2)(1-t)(1-t^2)} z^2 + \cdots,
	\]
	and the first terms of its plethystic logarithm read:
	\begin{equation}\label{eq:plog example}
		\pLog \Omega_k(z; q,t) = -\frac{1}{(1-q)(1-t)}\left((-1)^{k+1} z + 
		z^2 \sum_{\substack{i\geq 0,\; j\geq 0,\; i+j\leq k-2, \\ i+j\equiv k\pmod{2}}} q^i t^j + \cdots\right).
	\end{equation}
	
	In \eqref{eq:omega example} one can recognize the localization formula for the bundle $O(k)$ on the Hilbert scheme (see \cite{haiman2002vanishing}). We have
	\[
	\Omega_k(z; q,t) = \sum_{n=0}^\infty (-1)^{kn} z^n \chi(\mathrm{Hilb}_n(\mathbb{C}^2), O(k)).
	\]
	
	Let $R_x=\mathbb{C}[x_1,\ldots,x_n]$, $R_{x,y}=\mathbb{C}[x_1,\ldots,x_n,y_1,\ldots,y_n]$. The symmetric group $S_n$ acts on these, permuting $x_i$ and $y_i$ simultaneously. For an $S_n$-module $M$ we denote by $M^+$ and $M^-$ the subspace of invariant and anti-invariant elements respectively. Notice that $R_{x,y}^-\subset R_{x,y}$ is a module over $R_{x,y}^+$. $(R_{x,y}^-)^k$ denotes the linear span of $k$-fold products of elements of $R_{x,y}^-$ inside $R_{x,y}$. This is also a module over $R_{x,y}^+$. Haiman obtained the following results about the bundles $O(k)$:
	\begin{theorem}[\cite{haiman2002vanishing}]
		For any $k\geq 0$ we have an isomorphism of $R_{x,y}^+$-modules 
		\[
		\Gamma(\mathrm{Hilb}_n(\mathbb{C}^2), O(k)) \cong (R_{x,y}^-)^k \subset R_{x,y},
		\]
		and the $i$-th cohomology of $O(k)$ vanishes for $i>0$. $(R_{x,y}^-)^k$ is free as a module over $R_x^+$.
	\end{theorem}
	
	For instance, for $k=0$ and $k=1$ we have 
	\[
	\Gamma(\mathrm{Hilb}_n(\mathbb{C}^2), O) = R_{x,y}^+, \qquad
	\Gamma(\mathrm{Hilb}_n(\mathbb{C}^2), O(1)) = R_{x,y}^-.
	\]
	
	From this we can obtain the following presentation for the dual of the local cohomology:
	\begin{prop}
		Let  $M_{n,k}:=H^n_L(O(k))\;\check{}\;$. Let $\Delta_x := \prod_{i<j\leq n} (x_i - x_j)$. For $k\geq 0$ we have 
		\[
		M_{n,k} = \left\{f\in \frac{1}{\Delta_x^k} R_{x,y}^- \;\Big|\; D(x_1,\ldots,x_n,\partial_{y_1},\ldots,\partial_{y_n}) f\in R_{x,y} \; (\forall D\in (R_{x,y}^-)^k) \right\}.
		\]
		For $k\geq 1$ we have
		\[
		M_{n,k} = \left\{f\in \frac{1}{\Delta_x^{k-1}} R_{x,y}^+ \;\Big|\; D(x_1,\ldots,x_n,\partial_{y_1},\ldots,\partial_{y_n}) f\in R_{x,y} \; (\forall D\in (R_{x,y}^-)^{k-1}) \right\}.
		\]
	\end{prop}
	For instance, we have $M_{n,0}=R_{x,y}^-$ and $M_{n,1}=R_{x,y}^+$. Equivalence of the two presentations follows from the fact that if $f\in \frac{1}{\Delta_x} R_{x,y}^-$ satisfies $Df\in R_{x,y}$ for all $D\in R_{x,y}^-$, then we have $f\in R_{x,y}^+$.
	
	In the case $n=2$ the module $R_{x,y}^-$ is generated over $R_{x,y}^+$ by $x_1-x_2$ and $y_1-y_2$. So the condition on $f$ is equivalent to 
	\[
	(x_1-x_2)^i (\partial_{y_1} - \partial_{y_2})^j f \in R_{x,y} \qquad(i+j=k-1).
	\]
	Equivalently, in the expansion
	\[
	f = \sum_{i,j} (x_1-x_2)^i (y_1-y_2)^j f_{ij}(x_1+x_2, y_1+y_2)
	\]
	we have $f_{ij}=0$ unless $i\geq 0$ or $i\geq -k+1+j$. Note that also $f_{ij}=0$ unless $i+j\equiv k-1\pmod{2}$. So $M_{2,k}$ has a direct sum decomposition $M_{2,k} = R_{x,y}^{(-1)^{k-1}} \bigoplus M_{2,k}'$, where $M_{2,k'}$ has the Hilbert series 
	\[
	\frac{1}{(1-q)(1-t)}\sum_{\substack{i\geq 1,\,j\geq 0,\,i+j\leq k-1\\i+j\equiv k-1\pmod{2}}} q^{-i} t^j,
	\]
	which exactly matches the $z^2$-term of \eqref{eq:plog example} after the substitution $q\to q^{-1}$.
	
	More generally, we notice that there is a \emph{shuffle product} $M_{n',k} \otimes M_{n'',k} \to M_{n'+n'',k}$ given by
	\[
	f\ast g =  
	\]
	\[
	\sum_{\sigma\in S_{n'+n''}/(S_n'\times S_{n''})} (-1)^{(k-1)|\sigma|} \sigma \left(f(x_1,\ldots,x_{n'},y_1,\ldots,y_{n'}) g(x_{n'+1},\ldots,x_{n'+n''},y_{n'+1},\ldots,y_{n'+n''})\right).
	\]
	Supported by computer experiments, we formulate the following conjecture.
	\begin{conjecture}\label{conj:iso vector spaces}
		Denote by $M_{n,k}^\textrm{shuffles}\subset M_{n,k}$ the span of all shuffle products of the form $f\ast g$ for $f\in M_{n',k}$ and $g\in M_{n-n',k}$, $0<n'<n$. Let the space of primitive elements be the quotient $M_{n,k}^\textrm{prim}:=M_{n,k} / M_{n,k}^\textrm{shuffles}$. We have an isomorphism of bigraded vector spaces
		\[
		H^*(X)\left[\sum_{i=1}^n x_i, \sum_{i=1}^n y_i\right] \cong M_{n,k}^\textrm{prim}.
		\]
	\end{conjecture}
	
	The ring structure of the cohomology is not immediately evident from this description, but our COHA tells us where to look for it. 
	\begin{prop}
		The Lie algebra $\mathscr{H}_2 \otimes H^*(\BP^1)$ acts on $M_{n,k}$ by the following operations:
		\[
		(T_{r,s}\otimes \mathrm{pt}) f  = \sum_{i=1}^n x_i^{s} \partial_{y_i}^r f,
		\]
		\[
		(T_{r,s}\otimes 1) f  = \sum_{i=1}^n x_i^{s-1} \partial_{y_i}^{r-1} (-s \partial_{y_i} y_i - r x_i \partial_{x_i}) f.
		\]
		This action is compatible with the shuffle product and therefore descends to $M_{n,k}^\textrm{prim}$.
	\end{prop}
	
	Now we describe the ring structure. First we expect that the following is true:
	\begin{conjecture}
		The space $M_{n,k}^\textrm{prim}$ is spanned by products of $\frac{1}{\Delta_x^{k-1}}$ and polynomials of the form
		\[
		\psi_s(\mathrm{pt}) := \sum_{i=1}^n x_i^s, \qquad \psi_s(1) := s \sum_{i=1}^n x_i^{s-1} y_i.
		\]
	\end{conjecture}
	This conjecture implies that $M_{n,k}^\textrm{prim}$ is isomorphic to the quotient of the ring of polynomials in the generators $\psi_s(\mathrm{pt})$, $\psi_s(1)$ by a certain ideal. So it endows $M_{n,k}^\textrm{prim}$ with a ring structure and it is natural to formulate
	
	\begin{conjecture}\label{conj:iso rings}
		The isomorphism of Conjecture \ref{conj:iso vector spaces} is an isomorphism of rings.
	\end{conjecture}
	
	\subsection{Some remarks}
	
	The Hilbert scheme has been the source of interesting bigraded modules. For instance, there was $(n+1)^{n-1}$ conjecture and $n!$ conjecture proved by Haiman \cite{haiman2002vanishing}. The $(n+1)^{n-1}$ conjecture is about the module of diagonal coinvariants.
	\begin{theorem}\cite{haiman2002vanishing}
		\[
		\dim R_{x,y}/(R_{x,y}^+)_{>0} = (n+1)^{n-1}
		\]
	\end{theorem}
	An explicit combinatorial presentation of the bigraded Frobenius character of $R_{x,y}/(R_{x,y}^+)_{>0}$ was known as a shuffle conjecture \cite{haglund2005combinatoriala} and was proved in \cite{carlsson2018proof}. Also in this context Haiman observed that the vector space is in a certain sense generated by the differential operators $\sum_{i=1}^n x_i^{s-1} \frac{\partial}{\partial y_i}$, which led to the operator conjectures in \cite{haiman1994conjectures}. The double coinvariants has been interpreted in terms of Borel-Moore cohomology of affine Springer fibers in \cite{carlsson2018affine}.
	
	Our story suggests that there are many more interesting bigraded modules to be constructed using the Hilbert scheme.
	
	\section{Poisson geometry of Coulomb branches}
	We would like to conclude by explaining a curious formalism which hints at how Conjectures \ref{conj:iso vector spaces}--\ref{conj:iso rings} can be extended to surfaces of higher genus. First, observe that, in view of \eqref{eq:geometric engineering}, we should consider the bundle of multiforms, i.e. the $g$-fold tensor product of the bundle of differential forms on the Hilbert scheme. How does the Lie algebra $\mathscr{H}_2 \otimes H^*(C)$ act? Quite naturally, in fact. Take a function $f\in R_{x,y}^+$, for instance $f=\sum_{i=1}^n x_i^s y_i^r$. Symplectic geometry associates to $f$ the Hamiltonian vector field $H_f$. Let
	\[
	\psi_f(\mathrm{pt}) := \text{multiplication by $f$},
	\]
	\[
	\psi_f(1) := \text{Lie derivative along $H_f$}.
	\]
	This is what we had in genus $0$. But with multiforms, additional operations arise. For $i=1,\ldots,g$ set
	\[
	\psi_f(\sigma_i) = (d f)^{(i)}.
	\]
	This is the operation of multiplication by $df$ acting on the $i$-th component in the $g$-fold tensor product. Set
	\[
	\psi_f(\sigma_{i+g}) = \iota^{(i)}_{H_f},
	\]
	the operation of insertion of $H_f$ into the $i$-th component. From 
	\[
	\psi_{fg}(1) = \mathrm{Lie}_{f H_g + g H_f}, \qquad
	\mathrm{Lie}_{f H_g} = f \mathrm{Lie}_{H_g} + \sum_{i=1}^g (d f)^{(i)} \iota^{(i)}_{H_g}
	\]
	we obtain the following identity:
	\[
	\psi_{fg}(1) = \psi_f(1) \psi_g(\mathrm{pt}) + \psi_f(\mathrm{pt}) \psi_g(1) + \sum_{i=1}^{2g} \psi_f(\sigma_i) \psi_g(\sigma_i^*).
	\]
	This formula precisely matches the K\"unneth decomposition of the diagonal class in $H^2(C\times C)$, which confirms the consistency of the formalism. Okounkov observed that our formalism essentially coincides with the Clifford algebras considered in \cite{kazhdan2024l}.
	
	\section*{Acknowledgments}
	This research is supported by the Consolidator Grant No. 101001159 ``Refined invariants in combinatorics, low-dimensional topology and geometry of moduli spaces'' of the European Research Council.
	
	\bibliographystyle{amsalpha}
	\bibliography{../refs}
	
\end{document}